\documentclass[12pt,reqno]{amsart}

\usepackage{enumerate}
\usepackage{hyperref}
\usepackage[arrow,matrix]{xy}

\DeclareMathOperator{\ann}{ann}%
\DeclareMathOperator{\cd}{cd}%
\DeclareMathOperator{\cor}{cor}%
\DeclareMathOperator{\Fix}{Fix}%
\DeclareMathOperator{\Ker}{ker}%
\DeclareMathOperator{\rank}{rank}%
\DeclareMathOperator{\res}{res}

\newcommand{\EP}{\chi}
\newcommand{\F}{\mathbb{F}}%
\newcommand{\Fc}{\mathfrak{F}}%
\newcommand{\Fp}{\F_p}%
\newcommand{\Gal}{\text{\rm Gal}}%
\newcommand{\Ic}{\mathcal{I}}%
\newcommand{\N}{\mathbb{N}}%
\newcommand{\pc}{\text{\rm pc}}%
\newcommand{\sep}{\text{\rm sep}}%
\newcommand{\Z}{\mathbb{Z}}%

\begin{document}

\title[Galois Module Structure]{Galois Module Structure of Galois
Cohomology and Partial Euler-Poincar\'e Characteristics}

\author[Nicole Lemire]{Nicole Lemire}
\author[J\'{a}n Min\'{a}\v{c}]{J\'an Min\'a\v{c}}
\address{Department of Mathematics, Middlesex College, \ University
of Western Ontario, London, Ontario \ N6A 5B7 \ CANADA}
\email{nlemire@uwo.ca}
\email{minac@uwo.ca}

\author[John Swallow]{John Swallow}
\address{Department of Mathematics, Davidson College, Box 7046,
Davidson, North Carolina \ 28035-7046 \ USA}
\email{joswallow@davidson.edu}

\begin{abstract}
    Let $F$ be a field containing a primitive $p$th root of unity,
    and let $U$ be an open normal subgroup of index $p$ of the
    absolute Galois group $G_F$ of $F$. Using the Bloch-Kato
    Conjecture we determine the structure of the cohomology group
    $H^n(U,\Fp)$ as an $\Fp[G_F/U]$-module for all $n\in\mathbb{N}$.
    Previously this structure was known only for $n=1$, and until
    recently the structure even of $H^1(U,\Fp)$ was determined only
    for $F$ a local field, a case settled by Borevi\v{c} and Faddeev
    in the 1960s.  For the case when the maximal pro-$p$ quotient
    $T$ of $G_F$ is finitely generated, we apply these results to
    study the partial Euler-Poincar\'e characteristics of $\EP_n(N)$
    of open subgroups $N$ of $T$. We show in particular that the
    $n$th partial Euler-Poincar\'e characteristic $\EP_n(N)$ is
    determined by only $\EP_n(T)$ and the conorm in $H^n(T,\Fp)$.
\end{abstract}


\date{October 20, 2006}


\thanks{Nicole Lemire was supported in part by Natural Sciences
and Engineering Research Council of Canada grant R3276A01. J\'an
Min\'a\v{c} was supported in part by Natural Sciences and
Engineering Research Council of Canada grant R0370A01, by a
Distinguished Professorship for 2004--2005 at the University of
Western Ontario, by the Mathematical Sciences Research
Institute, Berkeley, and by the Institute for Advanced Study,
Princeton. John Swallow was supported in part by National
Security Agency grant MDA904-02-1-0061.}

\maketitle

\newtheorem{theorem}{Theorem}
\newtheorem{proposition}{Proposition}
\newtheorem{lemma}{Lemma}
\newtheorem*{acknowledgements*}{Acknowledgements}
\newtheorem{corollary}{Corollary}

\theoremstyle{definition}
\newtheorem*{remark*}{Remark}

\parskip=10pt plus 2pt minus 2pt

Let $F$ be a field containing a primitive $p$th root of unity
$\xi_p$. Let $G_F$ be the absolute Galois group of $F$, $U$ an open
normal subgroup of $G_F$ of index $p$, and $G=G_F/U$.

In the 1960s Z.~I.~Borevi\v{c} and D.~K.~Faddeev classified the
possible $G$-module structures of the first cohomology groups
$H^1(U, \Fp)$ in the case $F$ a local field \cite{Bo}. Recently this
result was extended from local fields to all fields $F$ as above
\cite{MS1}, and in this more general context the result has been
further developed and applied in \cite{MS2}, \cite{MS3},
\cite{MSS1}, and \cite{MSS4}. It is important in the study of Galois
cohomology to extend these results to all cohomology groups
$H^n(U,\Fp)$, $n\in\N$.

The celebrated 1982 paper of Merkurjev-Suslin \cite{MeSu} achieved
substantial progress in the investigation of Galois cohomology, and
in hindsight we know that these results were already sufficient to
extend the results of Borevi\v{c} and Faddeev to $H^2(U,\Fp)$ for
any field $F$ containing a primitive $p$th root of unity.  Because
of connections with Brauer groups, Galois embedding problems, and
Galois pro-$p$-groups, a precise knowledge of the structure
$H^2(U,\Fp)$ as an $\Fp[G]$-module has significant applications: see
\cite{LLMS2}, for example, for consequences characterizing Galois
Demu\v{s}kin pro-$p$-groups.

The main conjectures in contemporary Galois cohomology include the
Bloch-Kato Conjecture, the natural generalization to Milnor
$K$-theory of Hilbert~90, and generalization of the exact sequences
in Milnor $K$-theory contained in \cite{MeSu} from the case $n=2$ to
any $n \in \N$.  V.~Voevodsky's's partially published proof of the
Bloch-Kato Conjecture, based upon earlier work of A.~S.~Merkurjev,
M.~Rost, and A.~A.~Suslin, encompasses all of these. (See \cite{R},
\cite{Su}, \cite{Vo1} and \cite{Vo2}.) In this paper we show that
two results related to the Bloch-Kato Conjecture, contained in
\cite{Vo1} and \cite{Vo2} and quoted precisely in
section~\ref{se:bkktheory}, are sufficient to determine the
structure of all $\Fp[G]$-Galois modules $H^n(U,\Fp)$ for all $n \in
\N$. Moreover, this structure depends only on simple arithmetical
invariants attached to the field extension $E/F$, where $E$ is the
fixed field of $U$ in the separable closure $F_\sep$ of $F$. In
particular, we use Hilbert~90 for Milnor $K$-theory in an essential
way to prove that no cyclic $\Fp[G]$-module of dimension $j$ with $2
< j < p$ can occur as a summand of $H^n(U,\Fp)$. Because we rely
only on the main conjectures mentioned, our results could already
have been obtained conditionally in the early 1980s.

We have already used the results on the $\Fp[G]$-structure of
$H^n(U,\Fp)$ to obtain a generalization of Schreier's formula to
higher Galois cohomology groups \cite{LLMS1}, to obtain a new
characterization of Demu\v{s}kin groups \cite{LLMS2}, and to
construct fields with free or trivial Galois cohomology modules
\cite{LMS}.  Moreover, using ideas in this paper, we clarify when an
analogue of Hilbert 90 is valid for Milnor $k$-theory and Galois
cohomology \cite{LMSS}.  Our results here also lead to a
characterization of all possible decompositions of $H^n(U,F_p)$ into
a sum of indecomposable $\Fp[G]$-modules. In a special case, this
has been achieved \cite[Theorem 4]{LLMS2}.  An investigation of
the moduli space of all decompositions of $H^n(U,\Fp)$ is a work in
progress.

Here we use the results to produce remarkably simple formulas for the
$n$th partial Euler-Poincar\'e characteristic of $U$.  These formulas
provide a new tool for studying the structure of Sylow $p$-subgroups
of absolute Galois groups, about which little is currently known.
Recall Artin-Schreier's important observation that the only finite
subgroups of absolute Galois groups are the trivial subgroup and
$\Z/2\Z$, and Becker's generalization of this result to the case of
the Galois group of a maximal $p$-extension \cite{Be}.  These results
lead naturally to the search for a classification of finitely
generated subgroups of Sylow $p$-subgroups of absolute Galois groups.
The question is settled for algebraic number fields and for fields of
transcendence degree $1$ over local fields \cite{E1}, \cite{E2},
\cite{E3}. In the case of $p=2$ for arbitrary fields $F$
with $\vert F^\times/F^{\times 2}\vert\le 8$ \cite{JW1} partial
progress was made, but there are still open questions about
possible Galois pro-$2$ groups with three generators.
Thus, we know only very few types of such groups, and while
the elementary type conjecture asserts that there are no others, there
is not much evidence for the validity of the conjecture.  (See
section~\ref{se:fpdim} for additional comments and references on this
conjecture.)  It is significant, therefore, that we are able to derive
results on the good behavior of partial Euler-Poincar\'e
characteristics---and its relation to cohomological
dimension---without appealing to the conjecture.

In the Appendix~\cite{BeLMS} written with Dave Benson, we use these
results to provide examples of pro-$p$-groups which cannot be
realized as absolute Galois groups.

\section{Main Theorems}

The main ingredient for our determination of the $G$-module
structure of $H^n(U,\Fp)$ is Milnor $K$-theory. (See \cite{Mi} and
\cite[Chap.~IX]{FV}.)  For $i\ge 0$, let $K_iF$ denote the $i$th
Milnor $K$-group of the field $F$, with standard generators denoted
by $\{f_1,\dots,f_i\}$, $f_1, \dots, f_i\in F\setminus \{0\}$. For
$\alpha\in K_iF$, we denote by $\bar\alpha$ the class of $\alpha$
modulo $p$, and we use the usual abbreviation $k_nF$ for
$K_nF/pK_nF$.  We denote by $(f)$ the element in $H^1(U,\Fp)$
corresponding to $f$ in $F\setminus \{0\}$ and by $(f_1,...,f_n)$
the cup product $(f_1)\cdots(f_n)$ in $H^n(U,\Fp)$.

Let $E$ be the fixed field of $U$ in the separable closure $F_\sep$
of $F$. We write $N_{E/F}$ for the norm map $K_nE\to K_nF$, and we
use the same notation for the induced map modulo $p$. We denote by
$i_E$ the natural homomorphism in the reverse direction. We also
apply the same notation $N_{E/F}$ and $i_E$ for the corresponding
homomorphisms between cohomology groups, although when convenient we
replace these with the equivalent maps $\cor_{E/F}$ and
$\res_{F,E}$. The image of an element $\alpha\in K_i F$ in
$H^i(G_F,\Fp)$ we also denote by $\alpha$. We formulate our results
in terms of Galois cohomology for intended applications, but we use
Milnor $K$-theory in our proof.

Our decomposition depends on four arithmetic invariants
$\Upsilon_1$, $\Upsilon_2$, $y$, $z$, which we define as follows.
Fix $a\in F$ such that $E=F(\root{p}\of{a})$, and let $\sigma\in G$
satisfy $\root{p}\of{a}^{\sigma-1}=\xi_p$. First, for an element
$\bar\alpha$ of $k_i F$, let
\begin{equation*}
    \ann_{k_{n-1}F}\bar\alpha = \ann
    \left(k_{n-1}F \xrightarrow{\bar\alpha \cdot -}
    k_{n-1+i}F\right)
\end{equation*}
denote the annihilator of the product with $\bar\alpha$. When the
domain of $\bar\alpha$ is clear, we omit the subscript on the map
and write simply $\ann\bar\alpha$.  We use analogous notation for
the annihilator of $\alpha$ in $H^i(U,\Fp)$ when we work in Galois
cohomology rather than in Milnor $K$-theory. Because we will often
use the elements $\overline{\{a\}}$, $\overline{\{\xi_p\}}$,
$\overline{\{a,a\}}$, and $\overline{\{a, \xi_p\}}$, we omit the
bars for these elements. We also omit the bar in the element
$\overline{\{\root{p}\of{a}\}}\in k_n E$.

Fix $n\in \N$ and $U$ an open normal subgroup of $G_F$ of index $p$,
and write $E = \Fix(U)$.  Define invariants associated to $E/F$ and
$n$ as follows:
\begin{align*}
    d &:= \dim_{\Fp}\ k_n F/N_{E/F}k_n E, \qquad\qquad\quad\ e :=
    \dim_{\Fp} N_{E/F}k_n E \\
    \Upsilon_1 &:= \dim_{\Fp} \ann_{k_{n-1}F}\{a, \xi_p\}/\ann
    \{a\}, \quad \Upsilon_2 := \dim_{\Fp} k_{n-1}F/ \ann\{a,\xi_p\}
\end{align*}
\begin{equation*}
    y :=
    \begin{cases}
        \dim_{\Fp}\ (N_{E/F}k_n E)\ /\ \{a\}\cdot k_{n-1}F, & p>2
        \\
        \dim_{\F_2}\ (N_{E/F}k_n E)\ /\ \{a\} \cdot
        \ann_{k_{n-1}F}\{a,-1\}, &p=2
    \end{cases}
\end{equation*}
\begin{equation*}
    z :=
    \begin{cases}
        \dim_{\Fp}\ (k_n F)\ /\ \left(\{\xi_p\}\cdot k_{n-1}F +
        N_{E/F}k_n E\right), &p>2 \\
        \dim_{\F_2}\ (k_n F)\ /\ \left(\{a\} \cdot k_{n-1}F +
        N_{E/F}k_n E\right), &p=2.
    \end{cases}
\end{equation*}

Our main results are then the following.

\begin{theorem}\label{th:pnot2}
    If $p>2$ and $n\in\N$ then
    \begin{equation*}
        H^n(U, \Fp) \simeq X_1 \oplus X_2 \oplus Y \oplus Z
    \end{equation*}
    where
    \begin{enumerate}
        \item\label{it:pn21} $X_1$ is a trivial $\Fp[G]$-module of
        dimension $\Upsilon_1$ and
        \begin{equation*}
            X_1\cap \res_{F,E} H^n (G_F,\Fp) = \{0\}
        \end{equation*}
        \item $X_2$ is a direct sum of $\Upsilon_2$ cyclic
        $\Fp[G]$-modules of dimension 2
        \item $Y$ is a free $\Fp[G]$-module of rank $y$
        \item\label{it:pn24} $Z$ is a trivial $\Fp[G]$-module of
        dimension $z$ and
        \begin{equation*}
            Z \subset \res_{F,E} H^n (G_F,\Fp).
        \end{equation*}
    \end{enumerate}
    Further we have
    \begin{enumerate}
        \setcounter{enumi}{4}
        \item $Y^G = \res_{F,E} \cor_{E/F} H^n(U,\Fp)$
        \item $\cor_{E/F}\colon X_1\oplus X_2 \to (a)\cdot
        H^{n-1}(G_F,\Fp)$ is surjective
        \item\label{it:pnot2i7} $\Upsilon_1+\Upsilon_2+y=e$
        \item\label{it:pnot2i8} $\Upsilon_2+z=d$.
\end{enumerate}
\end{theorem}

\begin{theorem}\label{th:p2}
    If $p=2$ and $n\in\N$ then
    \begin{equation*}
        H^n(U, \F_2) \simeq X_1 \oplus Y \oplus Z
    \end{equation*}
    where
    \begin{enumerate}
        \item\label{it:p21} $X_1$ is a trivial $\F_2[G]$-module of
        dimension $\Upsilon_1$ and
        \begin{equation*}
            X_1 \cap \res_{F,E} H^n (G_F,\F_2) = \{0\}
        \end{equation*}
        \item $Y$ is a free $\F_2[G]$-module of rank $y$
        \item\label{it:p23}\label{it:p2i3} $Z$ is a trivial $
        \F_2[G]$-module of
        dimension $z$ and
        \begin{equation*}
            Z \subset \res_{F,E} H^n (G_F,\F_2).
        \end{equation*}
    \end{enumerate}
    Further we have
    \begin{enumerate}
        \setcounter{enumi}{3}
        \item\label{it:p2i4} $Y^G = \res_{F,E} \cor_{E/F}
        H^n(U,\F_2)$
        \item\label{it:p2i5} $\cor_{E/F}\colon X_1 \to (a)\cdot \ann
        (a, -1)$ is an isomorphism
        \item\label{it:p2i6} $\Upsilon_1+y=e$
        \item $\Upsilon_2+z=d$.
    \end{enumerate}
\end{theorem}

\begin{remark*}
    The case $n=1$ in the two theorems recovers the results of
    \cite{MS1}.
\end{remark*}

The isomorphism classes of $\Fp[G]$-modules $X_1 \oplus Z$, $X_2$,
and $Y$ are determined by $H^n(U,\Fp)$, but the decomposition of
$H^n(U,\Fp)$ into summands $X_1$, $X_2$, $Y$, and $Z$ is not
canonical. However, there is an equivalent, canonical version of our
results, contained in the following statements. These statements
describe $(k_nE)^G$ and the intersection of $(k_nE)^G$ with images
of $k_nE$ under successive powers of the augmentation ideal of
$\Fp[G]$.  Here the augmentation ideal is generated by $\sigma-1$,
but the statements are independent of the choice of a generator
$\sigma$ of $G$ as well as a primitive $p$th root of unity $\xi_p$.
Assume first $p>2$. (See Proposition~\ref{pr:sigmamin1} and
Lemma~\ref{le:fixedes}.)
\begin{enumerate}
    \item for each $i \ge 3$, \\ \\ $(\sigma-1)^{i-1} k_nE \cap
    (k_nE)^G = i_E N_{E/F} k_nE = (\sigma-1)^{p-1} k_nE.$ \\
    \item $(\sigma-1) k_nE \cap (k_nE)^G = i_E (\{\xi_p\}\cdot
    k_{n-1}F) + i_E N_{E/F} k_nE$. \\
    \item\label{it:stat3} $0\to\ann\{a\}\to k_{n-1}F
    \xrightarrow{\{a\}\cdot -} k_n F \xrightarrow{i_E} (k_n E)^G
    \xrightarrow{N_{E/F}} \\ \{a\}\cdot \ann \{a,\xi_p\}\to 0$.
\end{enumerate}
If $p=2$ the exact sequence (\ref{it:stat3}) is equivalent to
Theorem~\ref{th:p2}. This equivalent, canonical reformulation of
Theorems~\ref{th:pnot2} and \ref{th:p2} follows from their proofs in
sections~\ref{se:bkktheory}, \ref{se:notelem}, and
\ref{se:prfsthms12}.

Now we turn to an application of the results above to Galois
pro-$p$-groups. For a pro-$p$-group $T$ with finite cohomology
groups $H^i(T,\Fp)$ for $0\le i\le n$, the $n$th partial
Euler-Poincar\'e characteristic $\EP_n(T)$ is defined as
\begin{equation*}
    \EP_n(T) = \sum_{i=0}^n (-1)^i \dim_{\Fp} H^n(T,\Fp).
\end{equation*}
Generalizing a previous result of \v{S}afarevi\v{c}, Koch observed a
very interesting criterion for cohomological dimension $\cd (G)$
being $\le n$, using the Euler-Poincar\'e characteristic of open
subgroups of $T$. (See \cite[Thm.~5.5]{Ko}.) For a strengthening of
Koch's result see \cite{Sc}.

Suppose now that $T = G_F (p)$ is the maximal pro-$p$-quotient of an
absolute Galois group $G_F$ of a field $F$ containing a primitive
$p$th root of unity, and suppose $T$ has finite rank.  In this paper
we show that for $N$ an open subgroup of $T$ of index $p$, the
single condition that $\EP_n (N) = p \EP_n (T)$ is equivalent
to the condition that $\cor : H^n (N,\Fp) \to H^n (T,\Fp)$ is
surjective. Using this result, we strengthen Koch's criterion in our
case, showing that if $\EP_n(N)=p\EP_n(T)$ for all open
subgroups $N$ of index $p$, then $\cd(T)\le n$.

The essential part of section~\ref{se:fpdim} is a formula for $\EP_n
(N)$, where $N$ is an open subgroup of $T$ of index $p$.  This formula
depends only on $p$, $\EP_n(T)$, and the conorm in degree $n$. The
connection with module structure takes particular notice of the
maximal free submodules in the cohomology groups; we show that
$\EP_n(N)$ differs from $p\EP_n(T)$ by a certain multiple of the
degree $n$ conorm, and if $p>2$ then $\EP_n(N)$ is determined by the
maximal free submodules in the cohomology groups together with this
conorm.  To see the role of the free submodules more clearly when the
dimensions are finite, we set $F^i(N)$ to be a maximal free
$\Fp[T/N]$-submodule of $H^i (N,\Fp)$ and define the $n$th partial
free Euler-Poincar\'e characteristic $\EP^{\Fc}_n(N)$ by
\begin{equation*}
    \EP^{\Fc}_n(N) = \sum_{i=0}^n (-1)^i \dim_{\Fp}
    F^i(N).
\end{equation*}
(Observe that $\EP^{\Fc}_n(N)$ is independent of the choice of
$F^i(N)$ for each $i$.)  Our theorem is as
follows.

\begin{theorem}\label{th:main}
    Let $p$ be a prime and $F$ a field containing a primitive $p$th
    root of unity.  Suppose $T=\Gal (F(p)/F)$ is a pro-$p$-group of
    finite rank, $N$ is a subgroup of $T$ of index $p$, and $n\in
    \N$.
    \begin{enumerate}[(a)]
        \item\label{it:ma} We have
        \begin{equation*}
            p\EP_n(T)-\EP_n(N)=(-1)^n(p-1) \dim_{\Fp}
            \frac{H^n(T,\Fp)}{\cor H^n(N,\Fp)}.
        \end{equation*}
        \item\label{it:mb}
        Now additionally assume either that $p>2$ or that $p=2$ and
        $-1\in N_{E/F}(E)$, where $E=\Fix(N)$. Then
        \begin{equation*}
            \EP_n(N) = \EP^{\Fc}_n(N) + (-1)^n \dim_{\Fp}
            \frac{H^n(T,\Fp)}{\cor H^n(N,\Fp)}.
        \end{equation*}
    \end{enumerate}
\end{theorem}

Under the initial hypotheses of Theorem~\ref{th:main} we deduce the
following corollaries.

\begin{corollary}\label{co:1}
    The following are equivalent.
    \begin{enumerate}
        \item $\EP_n(N)=p\EP_n(T)$
        \item $\cor H^n(N,\Fp)= H^n(T,\Fp)$
    \end{enumerate}
    Under the hypothesis either that $p>2$ or that $p=2$ and
    $-1\in N_{E/F}(E)$, then these are additionally equivalent
    to
    \begin{enumerate}
        \setcounter{enumi}{2}
        \item $\EP_n(N)=\EP_n^{\Fc}(N)$.
    \end{enumerate}
\end{corollary}

Corollary~\ref{co:2} below shows in particular that if $\EP_n (N)
= p\EP_n (T)$ for all open subgroups $N$ of $T$ of index $p$,
then the same relation holds for all open subgroups $V$ of index $p$
of any given open subgroup $U$ of $T$.

\begin{corollary}\label{co:2}
    Let $n\in \N$.  The following are equivalent:
    \begin{enumerate}
        \item\label{it:co21} $\EP_n(N)=p\EP_n(T)$ for all open
        subgroups $N$ of $T$ of index $p$
        \item\label{it:co22} $\cd(T) \le n$
        \item\label{it:co23} $\EP_n(V)=p\EP(U)$ for all open
        subgroups $U$ of $T$ and all open subgroups $V$ of $U$ of
        index $p$ in $U$.
    \end{enumerate}
\end{corollary}

\section{Bloch-Kato and Milnor $K$-theory}\label{se:bkktheory}

Our proofs rely on the following two results in Voevodsky's work on
the Bloch-Kato Conjecture. Because we apply Voevodsky's results in
the case when the base field contains a primitive $p$th root of
unity we shall formulate Voevodsky's results restricted to this
case. The first is the Bloch-Kato Conjecture itself:

\begin{theorem}[{\cite[Lemma~6.11 and \S 7]{Vo1} and
\cite[\S 6 and Thm.~7.1]{Vo2}}] \label{th:bk}\
    \begin{enumerate}
        \item Let $F$ be a field containing a primitive $p$th root
        of unity and $m\in\N$.  Then the norm residue homomorphism
        \begin{equation*}
           k_m F\to H^m(G_F, \mu_p)
        \end{equation*}
        is an isomorphism.

        \item For any cyclic extension $E/F$ of degree $p$, the
        sequence
        \begin{equation*}
            K_mE \xrightarrow{\sigma-1} K_mE \xrightarrow{N_{E/F}}
            K_mF
        \end{equation*}
        is exact.
    \end{enumerate}
\end{theorem}

\noindent The second result establishes an exact sequence
connecting $k_mF$ and $k_mE$ for consecutive $m$.  (We translate
the statement of the original result to $K$-theory using the
previous theorem.)  In the following result $a$ is chosen to
satisfy $E=F(\root{p}\of{a})$.

\begin{theorem}[{\cite[Def.~5.1 and Prop.~5.2]{Vo1}}]\label{th:es}
    Let $F$ be a field containing a primitive $p$th root of unity
    with no extensions of degree prime to $p$.  Then for any cyclic
    extension $E/F$ of degree $p$ and $m\in \N$, the sequence
    \begin{equation*}
        k_{m-1}E \xrightarrow{N_{E/F}} k_{m-1}F \xrightarrow{\{a\}
        \cdot -} k_m F \xrightarrow{i_E} k_m E
    \end{equation*}
    is exact.
\end{theorem}

We observe that we may remove the hypothesis that the field $F$ has
no extensions of degree prime to $p$.  We give a proof of the
following theorem in section~\ref{se:esext}.

\begin{theorem}[{Modification of Theorem~\ref{th:es}}]
\label{th:esext}
    Let $F$ be a field containing a primitive $p$th root of unity.
    Then for any cyclic extension $E/F$ of degree $p$ and $m\in \N$
    the sequence
    \begin{equation*}
        k_{m-1}E \xrightarrow{N_{E/F}} k_{m-1}F \xrightarrow{\{a\}
        \cdot -} k_m F \xrightarrow{i_E} k_m E
    \end{equation*}
    is exact.
\end{theorem}

\section{Notation and Lemmas}\label{se:notelem}

We fix $n\in \N$ and the cyclic extension $E=F(\root{p}\of{a})$, and
we write $k_{n-1}F = \ann\{a\} \oplus V \oplus W$, where $\ann\{a,
\xi_p\} = \ann\{a\} \oplus V$. Observe that $\Upsilon_1=\dim_{\Fp}
V$ and $\Upsilon_2=\dim_{\Fp} W$. We denote by $i_E\colon K_nF\to
K_nE$ the map induced by the inclusion of $F$ in $E$. In what
follows we will frequently refer to the element $\root{p}\of{a}$,
and so we abbreviate it by $A$. We recall that if $p=2$ then
$\{a,\xi_p\}=\{a,-1\}=\{a,a\}\in k_2F$, while if $p>2$ then $\{a,a\}
= 0\in k_2 F$.  Finally, we will use the projection formula for
taking the norms of standard generators of $K_i F$ (see
\cite[p.~81]{FW}).

\begin{lemma}\label{le:VW}
    We have a vector space isomorphism
    \begin{equation*}
        V\oplus W \xrightarrow{\{a\}\cdot -} \{a\}\cdot
        k_{n-1}F
    \end{equation*}
    and, if $p>2$, the compositum of the maps $\{\xi_p\}\cdot-$
    and $i_E$
    \begin{equation*}
        W \xrightarrow{\{\xi_p\}\cdot - } \{\xi_p\} \cdot W
        \xrightarrow{i_E}
        i_E(\{\xi_p\}\cdot W)
    \end{equation*}
    is a vector space isomorphism as well.
\end{lemma}

\begin{proof}
    The first isomorphism follows from the fact that $V\oplus W$ is
    a complement in $k_{n-1}F$ of the kernel of multiplication by
    $\{a\}$. For the second, assume $p>2$. Suppose that $\bar w\in
    W$ and $\bar\alpha = \{\xi_p\}\cdot \bar w\in \ker i_E$. Then by
    Theorem~\ref{th:esext}, $\bar\alpha = \{a\}\cdot \bar c$ for
    $c\in K_{n-1}F$.  Since $\{a,a\}=0$ we see that $\{a\}\cdot
    \bar\alpha = 0$.  But then $\bar w \in \ann \{a,\xi_p\}$ and so
    $\bar w=0$.
\end{proof}

For $\gamma\in K_nE$, let $l(\gamma)$ denote the dimension of the
cyclic $\Fp[G]$-submodule $\langle \bar \gamma \rangle$ of
$k_n E$ generated by $\bar\gamma$. Then we have
\begin{equation*}
    (\sigma-1)^{l(\gamma)-1}\langle \bar\gamma \rangle = \langle
    \bar\gamma\rangle^G \neq 0 \text{\ \ \ and \ \ }
    (\sigma-1)^{l(\gamma)}\langle \bar\gamma \rangle=0.
\end{equation*}

We denote by $N$ the map $(\sigma-1)^{p-1}$ on $k_n E$. Because
$(\sigma-1)^{p-1}=1+\sigma+\dots+\sigma^{p-1}$ in $\Fp[G]$, we may
use $i_EN_{E/F}$ and $N$ interchangeably on $k_n E$.

\begin{lemma}\label{le:extend} Let $p>2$.
    Suppose $\gamma\in K_nE$ with $l=l(\gamma)\ge 2$. Then if $l\ge
    3$,
    \begin{equation*}
        (\sigma-1)^{l-1}\bar{\gamma}\in i_{E}N_{E/F}k_nE,
    \end{equation*}
    and if $l=2$
    \begin{equation*}
        (\sigma-1)\bar{\gamma}\in i_E(\{\xi_p\}\cdot k_{n-1}F)+
        i_EN_{E/F}k_nE.
    \end{equation*}
\end{lemma}

\begin{proof}
    If $l=p$, then $(\sigma-1)^{p-1}k_nE=i_{E}N_{E/F}k_nE$
    shows the result in this case.

    Suppose $l<p$. Then $\bar{\gamma}\in \ker(\sigma-1)^{p-1}$ and
    so $i_{E}N_{E/F}(\bar{\gamma})=0$.  By Theorem~\ref{th:esext},
    there exists $b\in K_{n-1}F$ such that $N_{E/F}\bar{\gamma} =
    \{a\}\cdot \bar{b}$. Equivalently $N_{E/F}\gamma=\{a\}\cdot b+
    pf$ for some $f\in K_nF$. By the projection formula~(see
    \cite[p.~81]{FW}),
    \begin{equation*}
        N_{E/F}\big(\gamma-\{A\}\cdot i_{E}(b)-i_E(f)\big)=0.
    \end{equation*}
    Then by Theorem~\ref{th:bk}, there exists $\omega\in
    K_nE$ such that
    \begin{equation*}
        (\sigma-1)\omega=\gamma-(\{A\}\cdot i_{E}(b)+
         i_{E}f)
    \end{equation*}
    and hence, since $l\ge 2$,
    \begin{equation*}
        (\sigma-1)^{l-1}\bar{\gamma}=(\sigma-1)^l\bar{\omega}+
        (\sigma-1)^{l-2}i_{E} (\{\xi_p\}\cdot b).
    \end{equation*}

    If $l(\gamma)\ge 3$ we deduce
    \begin{equation*}
        (\sigma-1)^{l-1}\bar{\gamma}=(\sigma-1)^l\bar\omega
    \end{equation*}
    where $l(\omega)=l+1$. Set $\bar{\gamma}_{l+1}=\omega$ and
    repeat the argument. We obtain $\gamma_k\in K_nE$ of
    lengths $l< k\le p$ with
    \begin{equation*}
        (\sigma-1)^{l-1}\bar{\gamma}=(\sigma-1)^{k-1}\bar{\gamma}_k.
    \end{equation*}
    Take $\bar\alpha=\bar{\gamma}_p$ to obtain
    \begin{equation*}
        (\sigma-1)^{l-1}\bar{\gamma}\in i_{E}N_{E/F}k_nE,
    \end{equation*}
    as required.

    If $l(\gamma)=2$ we have that
    \begin{equation*}
        (\sigma-1)\bar{\gamma} = (\sigma-1)^2\bar{\omega} +
        i_{E}(\{\xi_p\}\cdot \bar{b})
    \end{equation*}
    for some $\omega\in K_nE$ and some $b\in K_{n-1}F $. We see
    that $l(\omega)\le 3$. If $l(\omega)<3$ then
    \begin{equation*}
        (\sigma-1)\bar{\gamma}\in i_E(\{\xi_p\}\cdot k_{n-1}F),
    \end{equation*}
    while if $l(\omega)=3$ then by the previous case we see that
    $(\sigma-1)^2\bar\omega\in i_{E}N_{E/F}k_nE$, so the result
    holds in either case.
\end{proof}

\begin{proposition}\label{pr:sigmamin1}
    If $p>2$, then
    \begin{equation*}
        (\sigma-1)k_nE\cap (k_nE)^G = i_E(\{\xi_p\}\cdot k_{n-1}F) +
        i_{E}N_{E/F}k_nE
    \end{equation*}
    and for $3\le i\le p$,
    \begin{equation*}
        (\sigma-1)^{i-1}k_nE\cap (k_nE)^G=i_EN_{E/F}k_nE.
    \end{equation*}

    If $p=2$ then
    \begin{equation*}
        (\sigma-1)k_nE\cap (k_nE)^G=i_EN_{E/F}k_nE.
    \end{equation*}
\end{proposition}

\begin{proof}
    We have that for $2\le i\le p$,
    \begin{equation*}
        i_EN_{E/F}k_nE=(\sigma-1)^{p-1}k_nE\cap (k_nE)^G\subset
        (\sigma-1)^{i-1}k_nE\cap (k_nE)^G.
    \end{equation*}
    The statement for $p=2$ follows immediately. Now assume that
    $p>2$. Note also that
    \begin{equation*}
        i_E(\{\xi_p\}\cdot k_{n-1}F) =(\sigma-1)(\{A\}\cdot
        k_{n-1}F)\subset (\sigma-1)k_nE\cap (k_nE)^G.
    \end{equation*}
    Thus $i_E(\{\xi_p\}\cdot k_{n-1}F) + i_EN_{E/F}k_nE \in (\sigma
    - 1)k_nE \cap (k_nE)^G$ and the reverse inclusion follows from
    Lemma~\ref{le:extend}.

    Now we shall prove that if $3\le i\le p$, then
    \begin{equation*}
        (\sigma-1)^{i-1}k_nE \cap (k_nE)^G = i_EN_{E/F} k_nE.
    \end{equation*}
    Because $i_EN_{E/F}k_nE = (\sigma-1)^{p-1}k_nE$ we see that
    \begin{equation*}
        i_EN_{E/F}k_nE \subset (\sigma-1)^{i-1}k_nE \cap (k_nE)^G
    \end{equation*}
    for each $1\le i\le p$. In order to establish the reverse
    inclusion when $3\le i\le p$, consider $(\sigma-1)^{i-1}
    \bar{\gamma} \in (\sigma-1)^{i-1}k_nE\cap (k_nE)^G$.  We
    may assume without loss of generality that $l(\gamma)=i$, and
    then since $p>2$, we have $(\sigma-1)^{i-1}\bar\gamma
    \in i_E N_{E/F}k_nE$ by Lemma~\ref{le:extend}, as desired.
\end{proof}

In the following lemma we elongate the exact sequence of
Theorem~\ref{th:esext}.

\begin{lemma}\label{le:fixedes}
    The following sequence is exact:
    \begin{equation*}
        0\to\ann\{a\}\to k_{n-1}F \xrightarrow{\{a\}\cdot -} k_n F
        \xrightarrow{i_E} (k_n E)^G \xrightarrow{N_{E/F}} \{a\}\cdot
        \ann \{a,\xi_p\}\to 0.
    \end{equation*}
    Here the map $\ann\{a\}\to k_{n-1}F$ is the natural inclusion.
\end{lemma}

\begin{proof}
    We show first that $N_{E/F}((k_n E)^G)\subset \{a\} \cdot \ann
    \{a,\xi_p\}$.  Let $\bar \alpha\in (k_n E)^G$ and $\beta =
    N_{E/F}\alpha$. Since $i_E(N_{E/F}\bar\alpha) =
    (\sigma-1)^{p-1} \bar\alpha= 0$ we
    have that $\bar\beta = N_{E/F}\bar \alpha = \{a\}\cdot \bar b$
    for some $b\in K_{n-1}F$ by Theorem~\ref{th:esext}.

    Suppose $p=2$. Since $\bar\beta$ is in the image of $N_{E/F}$,
    we have by Theorem~\ref{th:esext} that $\{a\} \cdot \bar \beta
    = \{a,a\} \cdot \bar b = 0$.  Since $\{a,a\} = \{a,-1\}$, we
    have $\bar b \in \ann \{a,-1\}$.

    Now suppose that $p>2$.  Write $\beta = \{a\} \cdot b + pf$ for
    some $f\in K_nF$. Then by the projection formula (see
    \cite[p.~81]{FW}),
    \begin{equation*}
        N_{E/F}\big( \alpha - (\{A\} \cdot i_E(b) + i_E(f)) \big) =
        0.
    \end{equation*}

    By Theorem~\ref{th:bk}, there exists $\omega\in K_nE$ such
    that
    \begin{equation*}
        (\sigma-1)\omega = \alpha - (\{A\} \cdot i_E(b) - i_E(f))
    \end{equation*}
    and hence $(\sigma-1)^2\bar\omega = \{\xi_p\} \cdot
    i_E(\bar{b})$.

    If $(\sigma-1)^2\bar\omega = 0$ then by Theorem~\ref{th:esext},
    $\{\xi_p\} \cdot \bar b = \{a\} \cdot \bar h$ for some $h\in
    K_{n-1}F$.  Because $\{a,a\}=0$, the right-hand side of the
    preceding equation is annihilated by $\{a\}$.  Therefore $\bar
    b\in \ann \{a,\xi_p\}$.

    If $(\sigma-1)^2\bar \omega \neq 0$ then $l(\omega)=3$ and
    Lemma~\ref{le:extend} shows that
    \begin{equation*}
        i_E(\{\xi_p\}\cdot \bar b)  =
        i_E(N_{E/F}(\overline{\lambda}))
    \end{equation*}
    for some $\lambda\in K_nE$.  By Theorem~\ref{th:esext}, we have
    $\{\xi_p\} \cdot \bar b = N_{E/F}(\overline{\lambda}) + \{a\}
    \cdot \bar h$ for some $h\in K_{n-1}F$.  Now by
    Theorem~\ref{th:esext} and the fact that $\{a,a\}=0$, the
    right-hand side of the preceding equation is annihilated by
    $\{a\}$.  Then $\bar b\in \ann \{a,\xi_p\}$.   Hence in all
    cases $N_{E/F}\bar\alpha \in \{a\}\cdot \ann \{a, \xi_p\}$.

    Exactness at the first two terms is obvious, and exactness at
    the third term follows from Theorem~\ref{th:esext}.

    For exactness at the fourth term, suppose $$\bar \gamma\in (k_n
    E)^G\mbox{ and }N_{E/F}\bar\gamma = 0.$$ Then $N_{E/F}\gamma=
    pf$ for $f\in K_nF$.  Let $\beta=\gamma-i_E(f)$.  Then
    $N_{E/F}\beta=0$ and by Theorem~\ref{th:bk} there exists
    $\alpha\in K_nE$ such that $(\sigma-1)\alpha = \beta$. If $p=2$
    then $\bar\beta=i_E(N_{E/F}\bar\alpha)\in i_E k_n F$ and we are
    done. Thus assume $p>2$.

    Now if $(\sigma-1)\bar\alpha=\bar\beta=0$ we are done as then
    $\bar\gamma=i_E(\bar f)$. Hence assume $(\sigma-1)\bar\alpha\neq
    0$. Then $l(\alpha)=2$ and by Lemma~\ref{le:extend} we see that
    \begin{equation*}
        \bar\beta=(\sigma-1)\bar\alpha\in\{\xi_p\}\cdot i_Ek_{n-1}F+
    i_EN_{E/F}k_nF\subset i_Ek_nF
    \end{equation*}
    and exactness at the fourth term is established.

    Finally we show the exactness at the fifth term. Since
    \begin{equation*}
        \{a\}\cdot\ann\{a,\xi_p\}=\{a\}\cdot V
    \end{equation*}
    it is enough to show that each element $\{a\}\cdot\bar{v}$ where
    $\bar{v}\in V$ can be written as $N_{E/F}\bar\alpha$ for some
    $\bar\alpha\in(k_n E)^G$. Observe that
    \begin{equation*}
        (\sigma-1)(\{A\}\cdot i_E\bar{v})=\{\xi_p\}\cdot
        i_E(\bar{v}).
    \end{equation*}
    Also we have
    \begin{equation*}
        N_{E/F}(\{A\}\cdot i_E (\bar v))=
        \begin{cases}
            \{a\}\cdot\bar v &\mbox{ if }p>2 \\
            \{-a\}\cdot\bar v &\mbox{ if }p=2.
        \end{cases}
    \end{equation*}
    Therefore it is enough to show that there exists an element
    $\bar\gamma\in k_n E$ such that $(\sigma-1)\bar\gamma =
    \{\xi_p\}\cdot i_E(\bar{v})$ and
    \begin{equation*}
        N_{E/F}\bar\gamma =
        \begin{cases}
            0 &\mbox{ if }p>2 \\
            \{-1\}\cdot\bar v &\mbox{ if }p=2.
        \end{cases}
    \end{equation*}
    Indeed then we can set $\bar\alpha=\{A\}\cdot
    i_E(\bar{v})-\bar\gamma$.

    Because $\bar{v}\in\ann\{a,\xi_p\}$ we see that $\{\xi_p\}\cdot
    i_E(\bar{v})\in\ann\{a\}$. By Theorem~\ref{th:esext} there
    exists $\bar\beta\in k_n E$ such that
    \begin{equation*}
        \{\xi_p\}\cdot\bar{v}=N_{E/F}\bar\beta\text{ and }
        i_E(N_{E/F}\bar\beta)=(\sigma-1)^{p-1}\bar\beta.
    \end{equation*}
    Then setting $\bar\gamma=(\sigma-1)^{p-2}\bar\beta$ we obtain
    our required element. The proof of our lemma is now complete.
\end{proof}

Next, we need a general lemma about $\Fp[G]$-modules.  The
straightforward proof of this lemma is omitted.

\begin{lemma}\label{le:excl}
    Let $M_i$, $i \in \Ic$, be a family of $\Fp[G]$-modules
    contained in a common $\Fp[G]$-module $N$.  Suppose that the
    $\Fp$-vector subspace $R$ of $N$ generated by all $M_i^G$ has
    the form $R = \oplus_{i \in \Ic} M_i^G$.  Then the
    $\Fp[G]$-module $M$ generated by $M_i$, $i \in \Ic$, has the
    form $M = \oplus_{i \in \Ic} M_i$.
\end{lemma}

Finally, we need a general structure proposition about
$\Fp[G]$-modules that shows that the structure of an $\Fp[G]$-module
$X$ can be determined from the structure of $X^G$.

\begin{proposition}\label{pr:fpgstruct}
    Let $X$ be an $\Fp[G]$-module.  Set $L_p=(\sigma-1)^{p-1}X$ and
    for $1\le i<p$, suppose that $L_i$ is an $\Fp$-complement of
    $(\sigma-1)^{i}X\cap X^G$ in $(\sigma-1)^{i-1}X\cap X^G$.

    Then there exist $\Fp[G]$-modules $X_i$, $i=1,\dots,p$, such
    that
    \begin{enumerate}
        \item $X=\bigoplus_{i=1}^{p} X_i$,
        \item $X_i^G=L_i$ for $i=1,\dots,p$,
        \item each $X_i$ is a direct sum of $\dim_{\Fp}(L_i)$ cyclic
        $\Fp[G]$-modules of length $i$, and
        \item for each $i=1, \dots, p$, there exists an
        $\Fp$-submodule $Y_i$ of $X_i$ with $\dim_{\Fp}(Y_i)=
        \dim_{\Fp}(L_i)$ such that $\Fp[G]Y_i=X_i$.
    \end{enumerate}
\end{proposition}

\begin{proof}
    Since $(\sigma-1)^{p-1}X\subset X^G$, we see that $L_p =
    (\sigma-1)^{p-1}X\cap X^G$. By reverse induction on $i=p, p-1,
    \dots,1$, we obtain
    \begin{equation*}
        (\sigma-1)^{i-1}X\cap X^G=\oplus_{j=i}^pL_j.
    \end{equation*}
    Observe that because $L_i$ is an $\Fp$-complement of
    $(\sigma-1)^{i}X\cap X^G$ in $(\sigma-1)^{i-1}X\cap X^G$, for
    each $i=1, \dots, p$, there exists an $\Fp$-submodule $Y_i$ of
    $X$ such that $(\sigma-1)^i Y_i = 0$ and $(\sigma-1)^{i-1}: Y_i
    \to L_i$ is an $\Fp$-isomorphism. Set $X_i$ to be the
    $\Fp[G]$-submodule of $X$ generated by $Y_i$.  By Exclusion
    Lemma \ref{le:excl}, the sum of the $X_k$ is direct since the
    sum of the $X_k^G=L_k$ is direct.

    We need only show that $X=\oplus_{k=1}^{p}X_k$. Let $J_k = \Ker
    (\sigma-1)^k \subset X$, $k=1,\dots,p$. Then we have the
    filtration $X^G = J_1 \subset J_2 \subset \dots \subset J_{p} =
    X$. Let $x\in X$. Then $(\sigma-1)^{p-1}x\in L_{p}$, and hence
    there exists $y_{p}\in Y_{p}$ with $(\sigma-1)^{p-1} x
    =(\sigma-1)^{p-1} y_{p}$. Therefore $x-y_{p}\in J_{p-1}$.

    More generally we show for all $k \in \{ 2,3,\dots,p \}$ that if
    $x \in J_k$ then there exist $x_i\in X_i$, $k\le i\le p$, such
    that $x-\sum_{i=k}^{p}x_i\in J_{k-1}$. We have already shown our
    statement for $k = p$. Hence assume that $x \in J_k$, $k < p$.
    Then
    \begin{equation*}
        (\sigma-1)^{k-1} x \in (\sigma-1)^{k-1}X\cap
        X^G=\oplus_{i=k}^ {p}L_i.
    \end{equation*}
    Hence there exist $l_i\in L_i$, $i=k,\dots,p$, such that
    \begin{equation*}
        (\sigma-1)^{k-1} x = \sum_{i=k}^{p}l_i,
    \end{equation*}
    and there exist $y_i\in Y_i$, $i=k,\dots,p$, such that
    $(\sigma-1)^{i-1}y_i = l_i$. Therefore
    \begin{equation*}
        (\sigma-1)^{k-1} x = (\sigma-1)^{k-1}\left(\sum_{i=k}^{p}
        (\sigma-1)^{i-k}y_i\right).
    \end{equation*}
    Setting $x_k=y_k\in X_k$, and $x_i = (\sigma-1)^{i-k}y_i\in
    X_i$, $i>k$, we have
    \begin{equation*}
        x-\sum_{i=k}^{p}x_i\in J_{k-1}
    \end{equation*}
    as required.

    Therefore, using the fact that
    \begin{equation*}
        J_1 = X^G = \oplus_{i=1}^{p} L_i = \oplus_{i=1}^{p}
        X_i^G \subset \oplus_{i=1}^{p} X_i,
    \end{equation*}
    as well as our statement above, we show by induction on $k$
    that $J_k \subset \oplus_{i=1}^{p} X_i$ for all $k =
    1,2,\dots,p$. In particular $J_{p} = X \subset
    \oplus_{i=1}^{p} X_i$, which completes our proof.
\end{proof}

\section{Proofs of Theorems~\ref{th:pnot2} and \ref{th:p2}}
\label{se:prfsthms12}

We determine now the structure of $k_nE$ as an $\Fp[G]$-module.  We
do so by invoking Proposition~\ref{pr:sigmamin1} to determine
$\Fp$-modules $L_i$, $i=1,\dots,p$, such that
\begin{equation*}
    L_p=(\sigma-1)^{p-1}k_nE=(\sigma-1)^{p-1}k_nE\cap (k_nE)^G,
\end{equation*}
and such that $L_i$ is an $\Fp$-complement of $(\sigma-1)^{i}
k_nE\cap (k_nE)^G$ in $(\sigma-1)^{i-1}k_nE\cap (k_nE)^G$.  In this
way we may then apply Proposition~\ref{pr:fpgstruct}.

\begin{proof}[Proof of Theorem~\ref{th:pnot2}]
    Assume that $p>2$, and set
    \begin{equation*}
        L_1 = X_1+Z, \quad
        L_2 = i_E(\{\xi_p\}\cdot W), \quad
        L_p = i_EN_{E/F}k_nE,
    \end{equation*}
    where we choose $\Fp$-complements $X_1$, $W$, and $Z$, as
    follows:
    \begin{align*}
        X_1 &\oplus i_Ek_nF &=&\ (k_nE)^G\\
        W &\oplus \ann\{a,\xi_p\} &=&\ k_{n-1}F \\
        Z &\oplus \left(i_E(\{\xi_p\} \cdot
        k_{n-1}F)+i_EN_{E/F}k_nE\right) &=&\ i_Ek_nF.
    \end{align*}
    Observe that $X_1 \cap i_Ek_nF = \{0\}$ and therefore $X_1$
    satisfies the second part of condition (\ref{it:pn21}) in
    Theorem~\ref{th:pnot2}.  Similarly, the second part of
    condition (\ref{it:pn24}) follows from the definition of $Z$.

    Further set $L_i=0$ if $i\neq 1$, $2$, $p$. We claim that the
    $\Fp$-modules $L_i$ satisfy the hypotheses of
    Proposition~\ref{pr:fpgstruct}.

    First observe that $L_p=i_EN_{E/F}k_nE=(\sigma-1)^{p-1}k_nE$
    satisfies the hypotheses of Proposition~\ref{pr:fpgstruct} for
    $L_p$.  Next, by Proposition~\ref{pr:sigmamin1}, if $i\ge 3$
    then
    \begin{equation*}
        (\sigma-1)^{i-1}k_nE\cap
        (k_nE)^G=i_EN_{E/F}k_nE=(\sigma-1)^{p-1}k_nE,
    \end{equation*}
    so that the $L_i=0$, $3\le i<p$, satisfy the hypotheses for
    $L_i$. Moreover, again by Proposition~\ref{pr:sigmamin1},
    \begin{equation*}
        (\sigma-1)k_nE\cap (k_nE)^G=i_E(\{\xi_p\}\cdot
        k_{n-1}F)+i_EN_{E/F}k_nE.
    \end{equation*}

    We show next that
    \begin{equation*}
        L_2\oplus L_p = i_E(\{\xi_p\}\cdot k_{n-1}F)+i_EN_{E/F}k_nE.
    \end{equation*}
    First let $\bar\gamma \in i_E(\{\xi_p\}\cdot k_{n-1}F)\cap
    i_EN_{E/F}k_nE$. Then $\bar{\gamma}=\{\xi_p\}\cdot i_E(\bar
    f)=i_EN_{E/F}\bar{\alpha}$ for some $f\in K_{n-1}F$ and some
    $\alpha\in K_nE$. By Theorem~\ref{th:esext},  $\{\xi_p\}\cdot
    \bar{f}-N_{E/F}\bar{\alpha}=\{a\}\cdot \bar b$ for some $b\in
    k_{n-1}F$. But then $\{a,\xi_p\}\cdot \bar{f}=0$, since
    $N_{E/F}k_nE=\ann\{a\}$ by Theorem~\ref{th:esext} and
    $\{a,a\}=0$. Therefore $\bar{f}\in \ann\{a,\xi_p\}$. Conversely,
    if $\bar f\in \ann \{\xi_p, a\}$ then $\{\xi_p\}\cdot \bar f \in
    N_{E/F}k_nE$. Because $W$ is an $\Fp$-complement of
    $\ann\{a,\xi_p\}$ in $k_{n-1}F$, we obtain that
    \begin{align*}
        i_E(\{\xi_p\}\cdot k_{n-1}F) +i_EN_{E/F}k_nE &=
        i_E(\{\xi_p\}\cdot W)\oplus i_EN_{E/F}k_nE \\ &=L_2\oplus
        L_p.
    \end{align*}

    Finally, by construction of $X_1$ and $Z$, observe that
    \begin{equation*}
        (k_nE)^G=X_1\oplus Z\oplus L_2\oplus L_p.
    \end{equation*}
    Therefore $L_1=X_1\oplus Z$ is an $\Fp$-complement of
    $(\sigma-1)k_nE\cap (k_nE)^G=L_2\oplus L_p$ in $(k_nE)^G$, as
    desired.

    Applying Proposition~\ref{pr:fpgstruct}, we obtain
    $\Fp[G]$-modules $X_2$ and $Y$ such that $(X_2)^G=L_2$,
    $Y^G=L_p=i_EN_{E/F}k_nE$, $X_2$ is a direct sum of
    $\dim_{\Fp}L_2$ cyclic $\Fp[G]$-modules of length 2, and $Y$ is
    a free $\Fp[G]$-module of rank $\dim_{\Fp}L_p$. By construction
    $X_1$ and $Z$ are trivial $\Fp[G]$-modules.

    We turn to the calculation of the dimensions of $X_1$, $Z$,
    $L_2$, and $L_p$.  Note that for $p>2$, $\ker(i_E)\subset
    N_{E/F}k_nE$. Indeed, by Theorem~\ref{th:esext},
    $\ker(i_E)=\{a\}\cdot k_{n-1}F$ and $N_{E/F}k_nE=\ann\{a\}$, and
    from $N_{E/F}\{A\} = \{a\}$ and the projection formula (see
    \cite[p.~81]{FW}) we obtain $\ker(i_E)\subset N_{E/F}k_nE$. Then
    we calculate
    \begin{align*}
        \dim_{\Fp}L_p &= \dim_{\Fp}i_EN_{E/F}k_nE
        =\dim_{\Fp}N_{E/F}k_nE/\ker(i_E) \\
        &=\dim_{\Fp}N_{E/F}k_nE/\{a\}\cdot k_{n-1}F \\ &=y; \\
        \dim_{\Fp}Z &= \dim_{\Fp}i_Ek_nF/i_E(\{\xi_p\}\cdot
        k_{n-1}F)+i_EN_{E/F}k_nE \\ &=\dim_{\Fp}k_nF/(\{\xi_p\}
        \cdot k_{n-1}F+N_{E/F}k_nE) \\ &=z; \qquad \text{and} \\
        \dim_{\Fp}L_2 &= \dim_{\Fp}i_E(\{\xi_p\}\cdot W)=\dim_{\Fp}W
        \\ &= \dim_{\Fp}k_{n-1}F/\ann\{a,\xi_p\}\\ &=\Upsilon_2.
    \end{align*}
    (In the determination of $\dim_{\Fp} L_2$ we use both
    Lemma~\ref{le:VW} and the definition of $W$.)

    As for $X_1$, by Lemma~\ref{le:fixedes} we have
    $\ker(N_{E/F})\cap (k_nE)^G=i_Ek_nF$ as well as
    $N_{E/F}((k_nE)^G)=\{a\}\cdot \ann\{a,\xi_p\}$. Then $X_1$, as
    an $\Fp$-complement of $i_Ek_nF$ in $(k_nE)^G$, is mapped
    isomorphically under $N_{E/F}$ onto $\{a\}\cdot
    \ann\{a,\xi_p\}=\{a\}\cdot V$. Therefore we calculate
    \begin{align*}
        \dim_{\Fp}X_1 &=\dim_{\Fp}\{a\}\cdot \ann\{a,\xi_p\}\\ &=
        \dim_{\Fp}\ann\{a,\xi_p\}/\ann\{a\}=\Upsilon_1.
    \end{align*}

    Next we claim that $N_{E/F}(X_1\oplus X_2)=\{a\}\cdot k_{n-1}F$.
    First observe that because $i_EN_{E/F}(X_1 \oplus X_2) =
    (\sigma-1)^{p-1}(X_1 \oplus X_2) = \{0\}$ we have
    \begin{equation*}
        N_{E/F}(X_1 \oplus X_2) \subset \ker i_E = \{a\}\cdot
        k_{n-1}F.
    \end{equation*}
    Since $\{a\}\cdot k_{n-1}F=\{a\}\cdot (V\oplus W)$, it suffices
    to show
    \begin{equation*}
        \{a\}\cdot V\subset N_{E/F}(X_1\oplus X_2) \text{ and }
        \{a\}\cdot W\subset N_{E/F}(X_1\oplus X_2).
    \end{equation*}
    By the above, $N_{E/F}$ maps $X_1$ isomorphically onto
    $\{a\}\cdot V$, so we have the first inclusion. For the second,
    let $\{a\}\cdot \bar{w}\in \{a\}\cdot W$. Then
    \begin{equation*}
        N_{E/F}(\{A\}\cdot i_E(\bar{w}))=\{a\}\cdot \bar{w}.
    \end{equation*}
    Observe
    that $X_2 = \Fp[G]Y_2$, where $Y_2$ is an $\Fp$-submodule of
    $X_2$ mapped isomorphically under $(\sigma-1)$ onto
    $i_E(\{\xi_p\}\cdot W)$. Therefore there exists $\bar{\gamma}\in
    Y_2\subset X_2$ such that
    \begin{equation*}
        (\sigma-1)\bar{\gamma} = (\sigma-1)(\{A\}\cdot i_E\bar{w}).
    \end{equation*}
    Then $\{A\}\cdot i_E\bar{w}-\bar{\gamma}\in (k_nE)^G$.  Now
    \begin{equation*}
        N_{E/F}((k_nE)^G)=\{a\}\cdot \ann\{a,\xi_p\}=\{a\}\cdot
        V=N_{E/F}(X_1).
    \end{equation*}
    Applying the projection formula, $\{a\}\cdot \bar{w} -
    N_{E/F}(\bar{\gamma})\in N_{E/F}(X_1)$.  Hence $\{a\}\cdot
    W\subset N_{E/F}(X_1\oplus X_2)$, as required.

    Lastly, we need to show the relations between the dimensions of
    these modules. From
    \begin{equation*}
        y=\dim_{\Fp}N_{E/F}k_nE/\{a\}\cdot k_{n-1}F
    \end{equation*}
    and
    \begin{align*}
        \dim_{\Fp}(\{a\}\cdot k_{n-1}F) &= \dim_{\Fp}\{a\}\cdot
        (V\oplus W)\\ & =\dim_{\Fp}(V\oplus W)=\Upsilon_1+\Upsilon_2
    \end{align*}
    (again using Lemma~\ref{le:VW}), we obtain that
    \begin{equation*}
        \Upsilon_1+\Upsilon_2+y=\dim_{\Fp}N_{E/F}k_nE=e.
    \end{equation*}
    Furthermore, since
    \begin{align*}
        z &= \dim_{\Fp}\ \frac{i_Ek_nF}{i_E(\{\xi_p\}\cdot
        k_{n-1}F)+i_EN_{E/F}k_nE} \\
        &=\dim_{\Fp}\ \frac{i_Ek_nF}{i_E(\{\xi_p\}\cdot W)\oplus
        i_EN_{E/F}k_nE}
    \end{align*}
    and $\Upsilon_2=\dim_{\Fp} i_E(\{\xi_p\}\cdot W)$, we use
    $\ker(i_E)\subset N_{E/F}k_nE$ to deduce
    \begin{equation*}
        \Upsilon_2+z=\dim_{\Fp}k_nF/N_{E/F}k_nE=d.
    \end{equation*}
    The proof is complete.
\end{proof}

\begin{proof}[Proof of Theorem~\ref{th:p2}]
    Assume $p=2$, and set
    \begin{align*}
        L_1 &= X_1+Z \\
        L_2 &= i_EN_{E/F}k_nE,
    \end{align*}
    where we choose $\F_2$-complements $X_1$ and $Z$, as follows:
    \begin{align*}
        X_1 &\oplus i_Ek_nF &=&\ (k_nE)^G \\
        Z &\oplus i_EN_{E/F}k_nE &=&\ i_Ek_nF.
    \end{align*}
    Observe that $X_1 \cap i_Ek_nE = \{0\}$ and therefore $X_1$
    satisfies the second part of condition (\ref{it:p21}) of
    Theorem~\ref{th:p2}. Similarly, the second part of condition
    (\ref{it:p23}) follows from the definition of $Z$.  Note also
    that by construction the sum $X_1+Z$ is direct.

    Since by Proposition~\ref{pr:sigmamin1}
    \begin{equation*}
        (\sigma-1)k_nE\cap (k_nE)^G=i_EN_{E/F}k_nE,
    \end{equation*}
    we see that $X_1+Z$ is an $\F_2$-complement of
    $(\sigma-1)k_nE\cap (k_nE)^G$ in $(k_nE)^G$. Therefore $L_1$ and
    $L_2$ satisfy the hypotheses of Proposition~\ref{pr:fpgstruct},
    and there exist a free $\F_2[G]$-module $Y$ with
    $Y^G=i_EN_{E/F}k_nE$ and trivial $\F_2[G]$-modules  $X_1$ and
    $Z$ such that
    \begin{equation*}
        k_nE=X_1\oplus Y\oplus Z.
    \end{equation*}

    We next determine the rank of $Y$ and the dimensions of
    $X_1$ and $Z$. We claim first that
    \begin{equation*}
        \ker(i_E)\cap N_{E/F}k_nE=\{a\}\cdot \ann\{a,-1\}.
    \end{equation*}
    Indeed, by Theorem~\ref{th:esext},  $\ker(i_E)=\{a\}\cdot
    k_{n-1}F$ and $N_{E/F}k_nE=\ann\{a\}$.  Then, we deduce from
    $\{a\}\cdot \bar f\in \ker(i_E)\cap N_{E/F}k_nE$ and the
    fact that $\{a,a\}=\{a,-1\}$ for $p=2$ that
    \begin{equation*}
        \{a,-1\}\cdot \bar f=\{a,a\}\cdot \bar f=0.
    \end{equation*}
    Conversely, suppose $\bar{f}\in \ann\{a,-1\}$.  Then because
    $\{a,a\}=\{a,-1\}$ we have
    \begin{equation*}
        \{a\}\cdot \bar{f}\in \{a\}\cdot k_{n-1}F\cap \ann\{a\}=
        \ker(i_E)\cap N_{E/F}k_nE.
    \end{equation*}

    We then calculate
    \begin{align*}
        \dim_{\F_2} L_2 &= \dim_{\F_2}i_EN_{E/F}k_nE \\ &=
        \dim_{\F_2} N_{E/F}k_nE/\{a\}\cdot \ann\{a,-1\}=y
    \end{align*}
    and
    \begin{align*}
        \dim_{\F_2} Z &= \dim_{\F_2}i_Ek_nF/i_EN_{E/F}k_nF \\ &=
        \dim_{\F_2} k_nF/(N_{E/F}k_nE+\{a\}\cdot k_{n-1}F) = z.
    \end{align*}

    As for $X_1$, by Lemma~\ref{le:fixedes}, we observe that
    $\ker(N_{E/F})\cap (k_nE)^G=i_Ek_nF$ and $N_{E/F}((k_nE)^G) =
    \{a\}\cdot \ann\{a,-1\}$. Hence $X_1$, as an $\F_2$-complement
    of $i_Ek_nF$ in $(k_nE)^G$, is mapped isomorphically under
    $N_{E/F}$ onto $\{a\}\cdot \ann\{a,-1\}=\{a\}\cdot V$. Therefore
    \begin{align*}
        \dim_{\F_2}X_1 &= \dim_{\F_2}\{a\}\cdot \ann\{a,-1\}\\
        &=\dim_{\F_2} \ann\{a,-1\}/\ann\{a\}=\Upsilon_1.
    \end{align*}

    Now we establish the relations between these dimensions. Since
    \begin{equation*}
        y=\dim_{\F_2}N_{E/F}k_nE/\{a\}\cdot \ann\{a,-1\},
    \end{equation*}
    we have
    \begin{equation*}
        \Upsilon_1+y=\dim_{\Fp}N_{E/F}k_nE=e.
    \end{equation*}
    Finally, we calculate
    \begin{align*}
        d &= \dim_{\F_2}k_nF/N_{E/F}k_nE \\
        &=\dim_{\F_2}k_nF/(N_{E/F}k_nE+\{a\}\cdot k_{n-1}F)+  \\
        &\qquad \dim_{\F_2}(N_{E/F}k_nE+\{a\}\cdot
        k_{n-1}F)/N_{E/F}k_nE \\ &=z+\dim_{\F_2}\{a\}\cdot
        k_{n-1}F/\{a\}\cdot \ann\{a,-1\} \\ &= z+\dim_{\F_2}
        k_{n-1}F/\ann\{a,-1\}=z+\Upsilon_2.
    \end{align*}
    The proof is complete.
\end{proof}

\section{Proof of Theorem~\ref{th:esext}}\label{se:esext}

For the case $p=2$ we have the long exact sequence of Galois
cohomology groups due to Arason \cite[Satz~4.5]{A}. Suppose then
that $p>2$, and assume first that $F$ is perfect. Let $S$ be any
Sylow $p$-subgroup of $G_F=\Gal(F_{\sep}/F)$, and set $L$ to be the
fixed field of $S$. Because $F$ is perfect, the separable closure
$F_{\sep}$ is identical to the algebraic closure $\bar F$, and hence
each finite extension of $L$ has degree a power of $p$. In
particular, all of the hypotheses of Theorem~\ref{th:es} are valid
for the field $L$ in place of $F$. Furthermore, $([L:F],p)=1$. (Here
we use basic properties of supernatural numbers and Sylow
$p$-subgroups. See \cite[Chapter~1]{Ser}.) Therefore if
$E=F(\root{p}\of{a})$ is a cyclic extension of $F$ of degree $p$, so
is $EL=L(\root{p}\of{a})$ over $L$.  By Theorem~\ref{th:es} we see
that the sequence
\begin{equation*}
    k_{m-1}EL\xrightarrow{N_{EL/L}} k_{m-1}L
    \xrightarrow{\{a\}\cdot-} k_m L\xrightarrow{i_{EL}} k_m EL
\end{equation*}
is exact for each $m\in\N$.

We claim that $i_L\colon k_mF \to k_mL$ is injective. Indeed,
suppose that $i_L(\alpha)=0$ for some $\alpha\in k_m F$. Then there
exists a finite subextension $M/F$ of $L/F$ such that
$i_M(\alpha)=0$. Then
\begin{equation*}
    0=N_{M/F}(i_M(\alpha))=[M:F]\alpha.
\end{equation*}
(See \cite[Thm.~XI.3.8]{FV}.) Because $[M:F]$ is coprime with $p$,
we see that $\alpha=0$ and $i_L$ is injective as asserted. Similarly
we have that $i_{EL}\colon k_m E\to k_m EL$ is injective.

We then have the following commutative diagram:
\begin{equation*}
\xymatrix{ k_{m-1}EL \ar[r]^{N_{EL/L}} & k_{m-1}L \ar[r]^{\{a\}\cdot
-} & k_m L \ar[r]^{i_{LE}} & k_m EL \\ k_{m-1}E \ar[r]^{N_{E/F}}
\ar[u]_{i_{EL}} & k_{m-1}F \ar[r]^{\{a\}\cdot -} \ar[u]_{i_{L}} &
k_m F \ar[r]^{i_E} \ar[u]_{i_{L}} &
k_m E \ar[u]_{i_{EL}} \\
}
\end{equation*}
The fact that the first square is commutative follows from
\cite[p.~383]{BT}. The commutativity of the remaining part of the
diagram is clear.  Because the vertical maps are injective, we see
that the bottom row of the diagram is a complex: the composition of
any two consecutive maps is the zero map. We now establish exactness
at the second and third terms of the complex.

Let $\alpha\in k_{m-1}F$ such that $\{a\}\cdot\alpha=0$. Then
$\{a\}\cdot i_L(\alpha)=0$ and therefore there exists an element
$\beta\in k_{m-1}EL$ such that $N_{EL/L}(\beta)=i_L(\alpha)$. Let
$M/F$ be a finite extension such that $\beta$ is defined over $EM$.
Then $N_{EM/M}(\beta)=i_M(\alpha)$ (see \cite[p.~383]{BT}), and we
have
\begin{equation*}
    N_{EM/F}(\beta) = N_{M/F}(N_{EM/M}(\beta)) = N_{M/F}(i_M
    (\alpha)) = [M:F]\alpha.
\end{equation*}
Thus $N_{E/F}(N_{EM/E}(\beta)) = N_{EM/F}(\beta) = [M:F]\alpha$.

Because $([M:F],p)=1$ we see that $\alpha\in N_{E/F}(k_{m-1}E)$.
Therefore we have established the exactness of our complex at
$k_{m-1}F$.

Now assume that $\alpha\in k_m F$ such that $i_E(\alpha)=0\in k_m
E$. Then arguing as above, we see that there exist a finite
extension $M/F$ and $\beta\in k_{m-1}M$ such that $\{a\}\cdot
\beta = i_M(\alpha)\in k_m M$.

Applying $N_{M/F}$ and using the projection formula we see that
\begin{equation*}
    \{a\}\cdot N_{M/F}(\beta) = N_{M/F}(i_M(\alpha)) =
    [M:F]\alpha.
\end{equation*}
Because $[M:F]$ is coprime with $p$, $\alpha\in\{a\}\cdot
N_{M/F}(\gamma)$ for a suitable element $\gamma\in k_{m-1}M$. Hence
we see that our complex is also exact at $k_m F$ and the full
complex is exact.

Now if $F$ is not perfect, let $F_{\pc}$ denote the perfect closure
of $F$ (see, for instance, \cite[pp.~69--70]{Ka}).  Since finite
purely inseparable extensions of $F$ are of $q$th-power degree,
where $q$ is the characteristic of $F$, and $q\neq p$ since
$\xi_p\in F$, we obtain $([F_{\pc}:F],p)=1$. The argument above
establishes the theorem for the perfect field $F_{\pc}$, and a
similar transfer argument descends from $F_{\pc}$ to $F$. \qed

\section{From $H^i(G_F,\Fp)$ to $H^i(G_F(p),\Fp)$}\label{se:from}

In the previous sections, we worked with cohomology groups of
absolute Galois groups.  In this section we observe that, assuming
the Bloch-Kato Conjecture, we may replace these groups with
cohomology groups of maximal pro-$p$-quotients of absolute Galois
groups.

Let $F(p)$ be the compositum of all finite Galois $p$-power
extensions of $F$ in a fixed separable closure $F_{\sep}$ of $F$,
$G_F = \Gal(F_{\sep}/F)$, and $T=\Gal(F(p)/F)=G_F(p)$. Since $F$
contains $\xi_p$ we see that $F(p)$ is closed under taking $p$th
roots and hence $H^1 (G_{F(p)},\Fp) = \{0\}$. By the Bloch-Kato
Conjecture (see Theorem~\ref{th:bk}) the subring of the cohomology
ring $H^* (G_{F(p)},\Fp)$ consisting of elements of positive degree
is generated by $H^1(G_{F(p)},\Fp)$. Hence for each $i \in \N,
H^i(G_{F(p)},\Fp) = \{0\}$. Then, considering the
Lyndon-Hochschild-Serre spectral sequence associated to the exact
sequence
\begin{equation*}
    1 \to G_{F(p)} \to G_F \to T \to 1,
\end{equation*}
we have that $\inf \colon H^i(T,\Fp) \to H^i(G_F,\Fp)$ is an
isomorphism for all $i \in \N \cup \{0\}$.

Recall that $E=F(\root{p}\of{a})$, and let $N=\Gal(F(p)/E)$.  By the
same argument, we have that $\inf \colon H^i(N,\Fp) \to
H^i(G_E,\Fp)$ is an isomorphism. A straightforward argument with
cochains shows that this isomorphism is $T/N$-equivariant.

We denote by $(a)\in H^1(T,\Fp)$ the element corresponding to $a$
and by $(a,\xi_p)$ the cup product $(a) \cdot (\xi_p)$ in
$H^2(T,\Fp)$.  From Theorems~\ref{th:bk}, \ref{th:esext}, and the
above argument we see that the following sequence is exact for $i
\ge 0$:
\begin{equation}\label{eq:sec5es}
    H^i(N,\Fp) \xrightarrow{\cor} H^i(T,\Fp) \xrightarrow{\
    -\cdot (a)} H^{i+1}(T,\Fp) \xrightarrow{\res}
    H^{i+1}(N,\Fp).
\end{equation}
Hence $\cor H^i(N,\Fp) = \ann_{H^i(T,\Fp)}(a)$. Observe that, by
naturality, the inflation map $\inf$ induces the following
isomorphisms:
\begin{equation*}
    \cor H^i(N,\Fp) = \ann_{H^i(T,\Fp)} (a) \stackrel{\inf}{\to}
    \cor H^i(G_E,\Fp) = \ann_{H^i(G_F,\Fp)} (a),
\end{equation*}
\begin{equation*}
    \ann_{H^i(T,\Fp)} (a,\xi_p) \stackrel{\inf}{\to} \ann_{H^i
    (G_F,\Fp)} (a,\xi_p).
\end{equation*}

From now on we shall freely use the isomorphisms related to the
cohomology groups $H^i(T,\Fp)$, $H^i(G_F,\Fp)$ and the cohomology
groups of their open subgroups observed in this section.

\section{Partial Euler-Poincar\'e Characteristics}\label{se:fpdim}

In this section we determine $\dim_{\Fp} H^i(N,\Fp)$ for all open
subgroups of $T$.  In the case when $T$ is finitely generated we use
this result to calculate the $n$th partial Euler-Poincar\'e
characteristic of $N$.  We show that a very simple relation exists
between the $n$th partial Euler-Poincar\'e characteristics of $N$ and
$T$, and we then sharpen Koch's well-known criterium for the
cohomological dimension of $T$.  These results help clarify the
structure of finitely generated subgroups of Sylow $p$-subgroups of
absolute Galois groups, a class about which we know comparatively
little.

Recall that the elementary type conjecture asserts a classification of
finitely generated subgroups of Sylow $p$-subgroups of absolute Galois
groups.  A version of this conjecture was first introduced in the
context of Witt rings of quadratic forms by Marshall \cite{M1} (see
also \cite{M2}), and Jacob and Ware then considered the conjecture in
terms of Galois groups $G_F(2)$ \cite{JW1}, \cite{JW2}.  Efrat
\cite{E1} extended the conjecture to possible groups $G_F(p)$: each
finitely generated group $G_F(p)$ may be constructed from groups
$\Z/2\Z$, $\Z_p$, and $G_F(p)$ for $F$ a local field, using some
specific semidirect products and free products in the category of
pro-$p$-groups.  (For a precise definition, see \cite{E1} and
\cite{E3}.  Another variant of the elementary type conjecture may be
found in \cite{Ku}. For a more recent investigation of the elementary
type conjecture, see \cite{En}.)

Although the elementary type conjecture attempts a far-reaching
generalization of the work of Artin-Schreier and Becker \cite{Be},
at the moment there is little evidence for the conjecture in such
generality, although the conjecture holds for algebraic number
fields and fields of transcendence degree $1$ over local fields
\cite{E1}, \cite{E2}, \cite{E3} and for $G_F(2)$ when $\vert
F^\times/F^{\times 2}\vert \le 8$ \cite{JW1}.  The significance of
this work on the $n$th partial Euler-Poincar\'e characteristic of
finitely generated subgroups and their connections with
cohomological dimension lies in the fact that we do not use the
elementary type conjecture to achieve the new results.  As a
consequence, the formulas provide a new tool for the investigation
of the validity of the elementary type conjecture. (For further
connections among the elementary type conjecture, Demu\v{s}kin
groups, and Galois modules, see \cite{LLMS2}.)

For the remainder of the section we will use without mention the
fact that $\cor H^i(N,\Fp) = \ann_{H^i(T,\Fp)} (a)$, which follows
from Theorem~\ref{th:esext}. Moreover, from this point on, we will
write $\ann (a)$ for $\ann_{H^i(T,\Fp)} (a)$ and $\ann (a,\xi_p)$
for $\ann_{H^i(T,\Fp)} (a,\xi_p)$.  For each $i\ge 0$ we set
\begin{align*}
    a_i &= \dim_{\Fp} \ann (a) \\
    d_i &= \dim_{\Fp} H^i(T,\Fp)/\ann (a)  \\
    h_i &= \dim_{\Fp} H^i(T,\Fp) = a_i+d_i.
\end{align*}

\begin{proposition}\label{pr:formal}
    Let $N$ be an open subgroup of $T = \Gal(F(p)/F)$ of index $p$. (We
    do not assume that $T$ has finite rank.)

    \begin{enumerate}[(a)]
        \item We have
        \begin{align*}
            \dim_{\Fp} H^i(N,\Fp) &= d_{i-1} + d_i + p \dim_{\Fp}
            \frac{\cor H^i(N,\Fp)}{(a)\cdot H^{i-1}(T,\Fp)}.
        \end{align*}
        \item Assume further either that $p>2$ or that $p=2$ and $a$
        is a sum of two squares in $F$.  Then
        \begin{align*}
            \dim_{\Fp} H^i(N,\Fp) &= d_{i-1} + d_i + py.
        \end{align*}
        \item Assume further either that $p>2$ and $a_i$ is finite,
        that $p=2$, $a$ is a sum of two squares in $F$, and $a_i$ is
        finite, or that $p=2$ and $h_{i-1}$ and $h_i$ are finite.  Then
        \begin{align*}
            \dim_{\Fp} H^i(N,\Fp) &= d_{i-1}+d_i+p(a_i-d_{i-1}).
        \end{align*}
    \end{enumerate}
\end{proposition}

\begin{proof}
    Assume first that $p>2$.  Observe that
    \begin{align*}
        d_{i-1} &= \dim_{\Fp} \frac{H^{i-1} (T,\Fp)}{\ann(a,\xi_p)}
        + \dim_{\Fp} \frac{\ann(a,\xi_p)}{\ann (a)} = \Upsilon_2 +
        \Upsilon_1.
    \end{align*}
    By Theorem~\ref{th:pnot2}(\ref{it:pnot2i8}), we have $d_i =
    \Upsilon_2+z$.  Then from Theorem~\ref{th:pnot2} we conclude
    \begin{align*}
        \dim_{\Fp} H^i(N,\Fp) &= \Upsilon_1 + 2\Upsilon_2 + py
        + z \\ &= (\Upsilon_1 + \Upsilon_2) + (\Upsilon_2 + z)
        +   p y \\ &= d_{i-1} + d_i + p y.
    \end{align*}
    In the case when $a_i< \infty$ we deduce
    \begin{equation*}
        \dim_{\Fp} H^i(N,\Fp) = d_{i-1} + d_i + p (a_i - d_{i-1}).
    \end{equation*}

    Now assume that $p=2$ and that $a$ is a sum of two squares in $F$.
    In this case, using Theorem~\ref{th:p2}, we have:
    \begin{equation*}
        \dim_{\F_2} H^i(N,\F_2) = \Upsilon_1 + 2y + z.
    \end{equation*}
    Because $a$ is a sum of two squares in $F$, we see that $-1 \in
    N_{E/F}(E^\times)$.  Then $(a)\cdot (a) = (a)\cdot (-1) = 0$ in
    $H^2 (T,\F_2)$ and we obtain that $\ann(a,-1) = H^{i-1}(T,\F_2)$
    and $\Upsilon_1 = \dim_{\F_2} H^{i-1}(T,\F_2)/\ann(a) = d_{i-1}$.
    Furthermore, in this case $(a)\cdot H^{i-1}(T,\F_2)\subset \cor
    H^i(N,\F_2)$ and therefore
    \begin{equation*}
        z = \dim_{\F_2} \frac{H^i(T,\F_2)}{\cor H^i(N,\F_2)} = d_i.
    \end{equation*}
    Thus we recover the same formula as above:
    \begin{equation*}
        \dim_{\F_2} H^i(N,\F_2) = d_{i-1} + d_{i} + 2y.
    \end{equation*}
    In the case when $a_i<\infty$ we deduce
    \begin{equation*}
      \dim_{\F_2} H^i (N,\F_2) = d_{i-1} + d_i + 2(a_i - d_{i-1}).
    \end{equation*}

    Finally assume that $p=2$ and that $a$ is not necessarily
    a sum of two squares in $F$. Then from Theorem~\ref{th:p2}
    \begin{align}
        \dim_{\F_2} H^i(N,\F_2) &= \Upsilon_1 + 2y + z \notag \\ &=
        \dim_{\F_2} \frac{\ann(a,-1)}{\ann(a)} + 2 \dim_{\F_2}
        \frac{\cor H^i(N,\F_2)}{(a)\cdot \ann(a,-1)} \notag \\
        &\qquad + \dim_{\F_2} \frac{H^i(T,\F_2)}{(a)\cdot
        H^{i-1}(T,\F_2) + \cor H^i(N,\F_2)}.\label{eq:prop3}
    \end{align}

    Now assume that $h_{i-1}$ and $h_i$ are finite.  Then
    equation~\eqref{eq:prop3} may be simplified, as follows.  Observe
    that
    \begin{align*}
      \dim_{\F_2} &\frac{(a) \cdot H^{i-1}(T,\F_2) + \cor H^i(N,\F_2)}
      {\cor H^i(N,\F_2)}  \\
        &= \dim_{\F_2} \frac{(a)\cdot H^{i-1}(T,
        \F_2)}{(a) \cdot H^{i-1}(T,\F_2) \cap \cor H^i (N,\F_2)}
    \end{align*}
    We claim that
    \begin{equation*}
      (a) \cdot H^{i-1}(T,\F_2) \cap \cor H^i (N,\F_2) =
      (a) \cdot \ann (a,-1).
    \end{equation*}
    The right-hand side is contained in the left-hand side by
    Theorem~\ref{th:p2}(\ref{it:p2i5}). Now suppose that for some
    $h\in H^{i-1}(T,\F_2)$, $(a)\cdot h=\cor s$ for some $s\in
    H^i(N,\F_2)$. Then $(a,-1)\cdot h=(a, a)\cdot h=0$ since $\cor
    H^i(N,\F_2)=\ann (a)$ and, since $p=2$, $(a,a)=(a,-1)$. Then
    $h\in \ann(a,-1)$ and our claim is proved.
    Therefore
    \begin{align*}
        \dim_{\F_2} \frac{(a) \cdot H^{i-1}(T,\F_2) + \cor
        H^i(N,\F_2)} {\cor H^i(N,\F_2)} &=
        \dim_{\F_2} \frac{(a) \cdot H^{i-1} (T,\F_2)}{(a)\cdot \ann
        (a,-1)}\\ &= \dim_{\F_2} \frac{H^{i-1}(T,\F_2)}{\ann
        (a,-1)}.
    \end{align*}

    Now in equation~\eqref{eq:prop3} we add and subtract two copies of
    \begin{equation*}
        \dim_{\F_2} \frac{(a) \cdot H^{i-1} (T,\F_2)}{(a)\cdot \ann
        (a,-1)}\\ = \dim_{\F_2} \frac{H^{i-1}(T,\F_2)}{\ann (a,-1)}.
    \end{equation*}
    We then obtain
    \begin{align*}
        \dim_{\F_2} H^i(N,\F_2) &= \left( \dim_{\F_2} \frac
        {\ann(a,-1)}{\ann(a)} + \dim_{\F_2} \frac{H^{i-1}(T,\F_2)}
        {\ann(a,-1)}\right) + \\ & \left( \dim_{\F_2}
        \frac{H^i(T,\F_2)} {\cor H^i(N,\F_2)} \right) \\ &+2
        \big(\dim_{\F_2} \cor H^i(N,\F_2) -\dim_{\F_2} (a)\cdot
        H^{i-1}(T,\F_2)\big) \\ &= d_{i-1} + d_i + 2 (a_i -
        d_{i-1}).
    \end{align*}
\end{proof}

\begin{remark*}
    We observe that if $p=2$ and $h_{i-1}$ and $h_i$ are finite, the
    last formula of the proposition also follows from Arason's long
    exact sequence \cite[Satz 4.5]{A}.
\end{remark*}

\section{Proofs of Theorem~\ref{th:main} and Corollary~\ref{co:2}}

\begin{proof}[Proof of Theorem~\ref{th:main}]
    Let $T=G_F(p)$ be of finite rank and $N$ an open subgroup of index
    $p$.  Hence the $n$th partial Euler-Poincar\'e characteristic
    $\EP_n(T)$ is well-defined.  For $0\le i\le n$, let $a_i$, $d_i$,
    and $h_i$ be defined as in section~\ref{se:fpdim}.  Because
    $\dim_{\Fp} F^\times/F^{\times p} < \infty$ implies that $\dim_{\Fp}
    E^\times/E^{\times p} < \infty$ by Theorems~\ref{th:pnot2} and
    \ref{th:p2} in the case $n=1$, we see that $a_i$, $d_i$, $h_i$,
    and $\EP_i (N)$ are all well-defined for all $i \in \N$.

    By Proposition~\ref{pr:formal} we calculate
    \begin{align*}
        \EP_n (N) &= \sum_{i=0}^n (-1)^i \dim_{\Fp} H^i(N,\Fp)\\
        &= 1 + \sum_{i=1}^n (-1)^i \big(d_{i-1}+d_i+p(a_i-d_{i-1})\big)
        \\ &= (-1)^n d_n + p\left(\sum_{i=1}^n (-1)^i a_i\right) -
        p\left(\sum_{i=1}^n (-1)^i d_{i-1}\right).
    \end{align*}
    Since $h_i=a_i+d_i$, we have
    \begin{align*}
      p\EP_n(T)-\EP_n(N) &= p\left(1+\sum_{i=1}^n (-1)^i(a_i+d_i)\right)
      + (-1)^{n+1}d_n \\ &\qquad - p\left(\sum_{i=1}^n (-1)^i a_i\right)
      + p\left(\sum_{i=1}^n (-1)^i d_{i-1}\right)\\
      &= p(-1)^nd_n+(-1)^{n+1}d_n \\
      &= (-1)^n(p-1)d_n.
    \end{align*}
    We have therefore shown item~\eqref{it:ma}.

    Now if $p>2$ then by Theorem~\ref{th:pnot2}(\ref{it:pnot2i7}) and
    the fact that $F^0(N)=\{0\}$ we have
    \begin{equation*}
      \chi^{\Fc}_n(N) = p \sum_{i=1}^{n} (-1)^i (a_i-d_{i-1}).
    \end{equation*}
    Moreover, if $p=2$ and $-1\in N_{E/F}(E)$, then $(a,-1)=0$ in
    $H^2(T,\F_2)$ and $\ann_{H^{n-1}(T,\F_2)} (a,-1) =
    H^{n-1}(T,\F_2)$.  By Theorem~\ref{th:p2}(\ref{it:p2i6}) we
    obtain that the formula for the rank of $Y$ has the same shape
    as the formula in the case $p>2$:
    \begin{align*}
        \rank_{\F_2[G]} Y &= \dim_{\F_2} \cor H^n(N,\F_2)/(a)\cdot
        H^{n-1}(T,\F_2) \\ &= a_n - \Upsilon_1 \\ &= a_n-d_{n-1}.
    \end{align*}
    Therefore in this case we have
    \begin{equation*}
        \chi^{\Fc}_n(N) = 2 \sum_{i=1}^{n} (-1)^i (a_i-d_{i-1})
    \end{equation*}
    as well.
    In either case $p>2$ or $p=2$ and $-1\in N_{E/F}(E)$, we then
    calculate
    \begin{align*}
        \EP_n (N) &= 1 + \sum_{i=1}^n (-1)^i \big(d_{i-1}+d_i+
        p(a_i-d_{i-1})\big) \\ &= (-1)^n d_n + \chi_n^\Fc(N),
    \end{align*}
    showing equation~\eqref{it:mb}.
\end{proof}

Corollary~\ref{co:1} is an immediate consequence of
Theorem~\ref{th:main}.

\begin{proof}[{Proof of Corollary~\ref{co:2}}]\

    (\ref{it:co21})$\implies$(\ref{it:co22}): Assume that
    $\EP_n(N)=p\EP_n(T)$ for all subgroups $N$ of $T$ of index
    $p$.  Then by Corollary~\ref{co:1}, the corestriction map
    $\cor_n:H^n(N,\Fp)\to H^n(T,\Fp)$ is surjective for all such
    subgroups $N$.  By exact sequence~\eqref{eq:sec5es} of
    section~\ref{se:from}, $\res_{n+1}:H^{n+1} (T,\Fp)\to
    H^{n+1}(N,\Fp)$ is injective for all such subgroups $N$.

    By the Bloch-Kato Conjecture (see Theorem~\ref{th:bk}),
    $H^{n+1}(T,\Fp)$ is generated by elements $(a_1) \cdot (a_2)
    \cdots (a_{n+1})$, $a_i\in F^\times$, $(a_i)\in H^1(T,\Fp)$. We
    claim that these elements are equal to zero.  If $(a_1)=0$, then
    $(a_1)\cdots (a_{n+1})=0$.  If $(a_1)\neq 0$, then let
    $E=F(\root{p}\of{a_1})$ be the cyclic extension of degree $p$
    and set $N=\Gal(F(p)/E)$.  Then $(a_1)\cdots (a_{n+1})$ lies in
    the kernel of $\res_{n+1}: H^{n+1}(T,\Fp)\to H^{n+1}(N,\Fp)$.
    Since $\res_{n+1}$ is injective, $(a_1) \cdots (a_{n+1}) = 0$.
    Thus $H^{n+1}(T,\Fp)=\{0\}$, as required.

    (\ref{it:co22})$\implies$(\ref{it:co23}) follows from the fact
    that for each open subgroup $U$ of $T$, $\cd(U) = \cd(T) \le n$
    (see \cite[Prop.~3.3.5]{NSW}) and hence for each open subgroup
    $V$ of $U$ of index $p$, $\cor\colon H^n(V,\Fp) \to H^n(U,\Fp)$
    is surjective. (See \cite[Prop.~3.3.8]{NSW}.) Then by
    Corollary~\ref{co:1} we have $\EP_n(V) = p\EP_n(U)$.

    (\ref{it:co23})$\implies$(\ref{it:co21}) follows by setting $U =
    T$.
\end{proof}

\section*{Acknowledgements}
We are grateful to M.~Rost and to J.-P.~Serre for correspondence
concerning the Bloch-Kato Conjecture and to C.~Weibel for his
counsel concerning the formulation of the main theorems.  We thank
A.~Adem, D.~J.~Benson, M.~Bush, J.~F.~Carlson, S.~Chebolu,
D.~Eisenbud, M.~Fried, F.~Hajir, H.~Koch, M.~Kolster, J.~Labute,
T.~Lawson, H.~Lenstra, C.~Maire, T.~Nguyen-Quang-Do, Z.~Reichstein,
M.~Rost, A.~Schultz, J.-P.~Serre, R.~T.~Sharifi, J.-P.~Tignol,
C.~Weibel, T.~Weigel, and P.~Zalesskii for their interesting
comments, suggestions, and encouraging remarks concerning this
paper. Finally we thank the referee for the very interesting and
useful comments and suggestions related to this paper. In particular
the Appendix~\cite{BeLMS} is a result of his/her suggestion.

\end{document}